# The uplift payment elimination through the Lagrangian relaxation of the redundant constraints

**Vadim BOROKHOV**
LLC "En+development"
**Russia**
vadimab@eurosib.ru

## SUMMARY

Integration of distributed energy sources (in particular, the renewable generation) to a power market often leads to an increase in volatility of the power output of the conventional power plants and a need to ensure additional reserves in the power systems. This, in turn, results in an introduction of the additional operating constraints on the conventional generation. Also, a large share of the zero marginal cost generators tend to lower the marginal prices for power. These factors make significant contributions to the market prices and the lost profit of the dispatchable generation that needs to be compensated to ensure the economic stability of the centralized dispatch outcome.

The centralized dispatch optimization problem of a power market with an integrated unit commitment procedure is generally non-convex. In certain cases, the non-convexities originate from the demand side as well. In the presence of non-convexities, the market may not have an equilibrium price for power that provides economic stability of the centralized dispatch outcome. In this case, to achieve an economically stable outcome, the uplift payments to the market players are introduced as part of the pricing principle. These payments compensate the market players for the lost profits – the profits foregone by accepting the centralized dispatch schedule. In the case of uniform (nodal) pricing for power, the total lost profit of the market players corresponds to the difference between the values of the original centralized dispatch optimization problem and the relaxed one. In particular, if the uniform (nodal) price for power is set using the convex hull pricing (CHP) method, the total lost profit of the market players equals the duality gap resulting from the Lagrangian relaxation of the power balance constraint.

Given the general pricing principle that involves the lost profit compensation in the form of the uplift payments, we study a question if it is possible to introduce new products/services at the market and the associated prices such that
1. the set of the optimal solutions to the centralized dispatch optimization problem is unaffected;
2. the profit received by each market player at the centralized dispatch schedule is unchanged;
3. no market player has a lost profit (i.e., the maximum profit of each market player on its private feasible set is attained at the centralized dispatch schedule).

Thus, no new lost profit emerges and the original lost profit of the market players is fully absorbed in the pricing for the additional products/services introduced at the market. The abovementioned products/services are described by the additively separable inequality constraints in the centralized dispatch optimization problem that undergo the Lagrangian relaxation procedure together with the power balance constraints. One of the possible ways to ensure condition #1 is to consider the redundant constraints, which are the constraints that hold on the feasible set of the original optimization problem. Therefore, the introduction of these constraints in the original constraint set does not affect the set of optimal solutions to the centralized dispatch optimization problem. Previous investigation has shown that introduction and subsequent Lagrangian relaxation of the affine (linear) redundant constraints do not change the duality gap.

This implies that if the market price is set using the CHP method, such constraints do not affect the aggregated maximum profit of the market players that could be attained on their private feasible sets. Since the redundant inequality constraints make only non-negative contributions to a market player's profit, the invariance of the aggregated maximum profit of the market players implies that the maximum profit of each market player is unchanged. For the case of the CHP method, prior research has illustrated that the Lagrangian relaxation of an affine redundant constraint may reduce the total uplift payment. (In such an approach, the total uplift payment no longer coincides with the duality gap.) However, in the general case, utilization of the affine constraints may not reduce the uplift payment to zero, while the introduction of the non-affine redundant constraints may change the duality gap (which means that the maximum values of the market player profits are modified).

We consider a special class of the (possibly non-linear) redundant constraints that are redundant not only on the feasible set of the centralized dispatch optimization problem (and, therefore, do not change the centralized dispatch outcome) but also on the larger set obtained by discarding the power balance constraint from the original constraint set. The Lagrangian relaxation of these redundant constraints may reduce the total uplift payment without changing the duality gap.

For any given market price (or a pricing algorithm that sets the producer revenue as a function of its output volume) in a uninode multi-period power market with fixed load, we explicitly construct a family of the redundant constraints that do not change the maximum profit of any producer and result in zero uplift payment. We show that the introduction and subsequent Lagrangian relaxation of just one redundant constraint in the centralized dispatch problem suffice to eliminate the uplift payments for all the producers.

Thus, the answer to the question above is positive, and a market player uplift payment can be fully absorbed by the redundant constraint relaxation. The analysis can be straightforwardly extended to a multi-node multi-period market with price-sensitive demand. We provide a number of examples (including a two-node power system with FTR-holders) illustrating an application of the proposed method.

## KEYWORDS

Uplift Reduction - Power pricing - Pricing non-convexities - Redundant constraint.

## 1. Introduction

Distributed energy production and energy storage solutions are the critical parts of the future power systems. Integration of the distributed renewable energy sources creates a number of challenges to the power system operation due to the intermittent nature of the solar and especially of the wind power production. Consequently, both net demand and a set of binding transmission constraints may differ at various stages of the power market planning such as the day–ahead market and the real-time market. In addition to development of efficient forecasting techniques for power output by these resources, this implies a need to introduce additional operating constraints on the dispatchable units in the systems with high penetration of renewables. In addition, physically realistic description of the power storage devices may involve non-convex constraints. These issues, combined with zero marginal cost of output by renewable energy sources, make a critical impact on the market prices (producing low or even negative market prices for power) and the associated economic mechanisms (e.g., the lost profit payments) utilized to stabilize the market outcome.

Many of deregulated electric power markets are based on a security-constrained centralized dispatch optimization problem with the objective function usually having the form of the total social welfare. This optimization problem also produces the system marginal price (or the locational marginal prices) for power [1]-[3]. If the centralized dispatch optimization problem is convex, then there is an equilibrium price for power (which is also a marginal price for power). In this case, all the market players acting as price-takers have no economic incentives to deviate from the centralized dispatch schedule. If the optimization problem is non-convex, then an equilibrium price may fail to exist (although, a marginal price may still be calculated) and no uniform price for power ensures the economic equilibrium of the centralized dispatch schedule. In this case, some alternative pricing approach should be implemented to ensure the economic stability of the centralized dispatch outcome [4]. The non-convexities may be present at the supply side (due to the no-load cost, start-up cost, non-zero minimum output limits, integral commitment decision variables, minimum up/down times, etc.) and from the flexible demand side as well [5]. Power markets in the USA address this problem by introducing the side-payments, while most of the power markets in Europe make no account for non-convexities in the power pricing mechanisms and the market players need to internalize them.

The power pricing problem in the presence of non-convexities has been extensively studied. In [6] it was proposed to introduce a new service (a unit being online) and the associated (unit-specific) prices by constraining the integral status variables to their values set by the centralized dispatch. However, these prices can be negative, and the resulting pricing method resembles the pay-as-bid approach. Bjørndal and Jörnsten in [7], [8] further improved this method by adding extra constraints and setting certain continuous variables to their optimal values as well. An alternative direction of research on the subject is given by the nonlinear (discriminatory) pricing methods for power with market player-specific prices [9]-[15]. Liberopoulos and Andrianesis [16] developed the minimum zero-sum uplift pricing approach that increases the price above the marginal costs and transfers all the additional payments (received by the profitable producers due to the price increase) to the producers with losses to make them whole. Araoz and Jörnsten [17] introduced a semi-Lagrangian relaxation approach that produces the same solution as the original centralized dispatch problem and results in a uniform market price with the non-confiscatory pricing for generators. Van Vyve [18] suggested a zero-sum uplift pricing scheme that minimizes the maximum contribution to the uplift financing in a market with the price-sensitive load. In the case of no price-sensitive load, this approach produces the market price that is equal to the maximum average cost of the producers. In the framework of a uniform market pricing, if an equilibrium price fails to exist, then to ensure the economic stability of the centralized dispatch outcome the pricing has to be supplemented with the lost profit (opportunity cost) compensations. In the context of the centralized treatment of non-convexities, these compensations are usually called the uplift payments without considering the other uplift components that are also present in the practical cases. The minimum-uplift pricing proposed in [19]-[21], also known as the convex hull pricing (CHP), produces a uniform market price that minimizes the total uplift payment, which is needed to ensure the economic stability of the centralized dispatch outcome. In this approach, for a given market price the lost profit of a market player is calculated as the difference between the maximum value of its profit function on the market player private feasible set and its profit received when following the centralized dispatch. In [22], [23] it was proposed to modify the minimum-uplift pricing method by excluding the power volumes that are not attainable in a decentralized market from the lost profit calculation since the opportunity to supply/consume these volumes is not forgone by a market player when accepting the centralized dispatch outcome. Since the private feasible sets used to calculate the lost profit are potentially reduced, this method produces the lower (or equal) total uplift payment compared to the CHP algorithm.

Eldridge et al. [24] obtained a bound on the redistributions that result from the sub-optimal schedules [25], [26] and showed that CHP minimizes the bound. Herrero et al. [27] showed that the uplift payments distort the uniform market pricing and decrease the transparency of the market; therefore, it is critical reducing these payments. Another direction of research is to modify the objective function of the centralized dispatch optimization problem to achieve lower total uplift. In this approach, the centralized dispatch outcome differs from that obtained from the social welfare optimization. To find the market prices that minimize the social welfare reduction due to schedules inconsistency and ensure non-negative generator profits a primal-dual approach was developed in [28]. In this method, however, some of the lost profit may not be compensated to producers, and the competitive equilibrium is not achieved. Huppmann and Siddiqui [29] proposed a general method to find economic equilibrium by modifying the centralized dispatch objective function to include the compensation payments. Fuller and Çelebi [30] proposed calculating the market price and the centralized dispatch outcome from the total opportunity cost minimization problem (with the power balance and power flow constraints included in the constraint set of the problem and the proper treatment of the self-dispatched demand). Shavandi et al. [31] extended the minimum complementarity approximation to the minimum total opportunity cost model to include price-responsive demand, reserve capacities, and reserve pricing and illustrated that the near-equilibrium solutions with little social welfare reduction may substantially lower the total opportunity cost.

For the CHP method, motivated by the fact that the Lagrangian relaxation of the linear redundant constraints does not change the duality gap [32], Zhang et al. [33] showed that introduction and the Lagrangian relaxation of one affine redundant constraint for each time instance reduce the total uplift payment, which compensates the lost profit of the market players. This new constraint is given by a linear combination of the power balance constraint and the capacity limit constraints of all the producers. The Lagrangian relaxation of this redundant constraint leads to the introduction of the new service (a unit being online) and the associated price in addition to the market price for power for each time instance. In this approach, the duality gap may no longer equal the total uplift, which is potentially reduced as some of the original uplift payment is absorbed into a payment received for being online. It is possible to show that in the general case the linear redundant constraints are not able to eliminate the total uplift payment. In [34] it was shown that if a redundant constraint is redundant not only on the feasible set of the centralized dispatch optimization problem but also on the feasible set of the relaxed problem (i.e., the feasible set defined by the original constraint set of the centralized dispatch optimization problem with the power balance constraint excluded), then the linearity requirement can be lifted. In this approach, the total uplift payment can be eliminated by the introduction of a non-linear redundant constraint in a general market pricing scheme.

In this paper, we study the problem of the total uplift payment (total lost profit) reduction in a general pricing setting, which fixes the producer revenue as some function of its status-output variables and the pricing parameters. Thus, our analysis is applicable to the cases with uniform pricing for power (such as marginal pricing, CHP) and non-linear pricing with the uplift payments. We extend the analysis of [34] to show its relation to the Lagrangian relaxation technique and dual optimization problem (if the price is set using the CHP method). For simplicity, we focus on a multi-period uninode power market with the fixed load. However, our analysis and the results can be straightforwardly translated to multi-node markets with price-sensitive demand. Without loss of generality, we suppose that each producer operates just one generating unit.

The paper is organized as follows. In Section 2, we introduce a set of redundant constraints and show how they can reduce the total uplift payment. We also observe that the same effect on the total uplift payment is achieved by just one redundant constraint. In Section 3, we formulate a stronger redundancy condition and present a general framework for the construction of a redundant constraint that reduces the total uplift to zero. The applications of the proposed approach are given in Section 4. The results are summarized in Section 5.

## 2. The total uplift reduction using the redundant constraints

Consider $T$-period uninode power market with the fixed load $\mathbf{d}$, $\mathbf{d} \in R_{\geq 0}^T$, with a set of producers $I$. For each producer $i \in I$ at time period $t \in \{1,..,T\}$, let $u_i^t$ and $g_i^t$ denote the commitment and output variables, respectively. Introduce $\mathbf{g}_i = (g_i^1,...,g_i^T)$, $\mathbf{u}_i = (u_i^1,...,u_i^T)$, $x_i^t = (u_i^t, g_i^t)$, $\mathbf{x}_i = (x_i^1,...,x_i^T)$, $\mathbf{x} = (\mathbf{x}_1,...,\mathbf{x}_{|I|})$, where $|\bullet|$ denotes cardinality of a set. Let $X_i$ be the producer $i$ private feasible set (which is assumed to be nonempty and compact) and $C_i(\mathbf{x}_i)$ be the offer cost function of the producer $i$. The centralized dispatch problem has the form

$$f^* = \min_{\substack{\mathbf{x}_i \in X_i, \forall i \in I \\ s.t. \sum_{i \in I} \mathbf{g}_i = \mathbf{d}}} \sum_{i \in I} C_i(\mathbf{x}_i). \quad (1)$$

Let $\Omega$ denote the feasible set of (1), which is assumed to be nonempty. Clearly, $\Omega$ is compact. Let $\mathbf{x}^* = (\mathbf{x}_1^*,...,\mathbf{x}_{|I|}^*)$ denote an optimal point of (1). If (1) has multiple optimal points, then $\mathbf{x}^*$ denote any one of them. (In what follows, we perform our analysis for a given choice of the optimal point.) We consider the case of uniform pricing for electric power, which is assumed to be the only commodity traded at the market. The Lagrangian relaxation of the power balance constraint $\sum_{i \in I} \mathbf{g}_i = \mathbf{d}$ gives the Lagrange function

$$L(\mathbf{p},\mathbf{x}) = \mathbf{p}^T(\mathbf{d} - \sum_{i \in I}\mathbf{g}_i) + \sum_{i \in I} C_i(\mathbf{x}_i) = \mathbf{p}^T\mathbf{d} - \sum_{i \in I}\pi_i^{st.}(\mathbf{p},\mathbf{x}_i)$$

with a multiplier $\mathbf{p} \in R^T$, which is interpreted as (the vector of) the uniform market prices, and the standard expression for a producer $i$ profit function $\pi_i^{st.}(\mathbf{p},\mathbf{x}_i) = R_i^{st.}(\mathbf{p},\mathbf{x}_i) - C_i(\mathbf{x}_i)$ with the revenue function $R_i^{st.}(\mathbf{p},\mathbf{x}_i) = \mathbf{p}^T\mathbf{g}_i$. Let us define $\pi_i^{st.*}(\mathbf{p}) = R_i^{st.}(\mathbf{p},\mathbf{x}_i^*) - C_i(\mathbf{x}_i^*)$ and $\pi_i^{st.+}(\mathbf{p}) = \max_{\mathbf{x}_i \in X_i} \pi_i^{st.}(\mathbf{p},\mathbf{x}_i)$, which are the producer $i$ profit at the centralized dispatch solution and its maximum profit attainable on the private feasible set $X_i$, respectively. Clearly, $\pi_i^{st.+}(\mathbf{p}) \geq \pi_i^{st.*}(\mathbf{p})$, $\forall i \in I$, and the difference $[\pi_i^{st.+}(\mathbf{p}) - \pi_i^{st.*}(\mathbf{p})]$ equals the profit lost by the producer $i$ when complying with the centralized dispatch solution. We assume that all the producers recover their lost profits through the uplift payments. For the total lost profit we have

$$\sum_{i \in I}[\pi_i^{st.+}(\mathbf{p}) - \pi_i^{st.*}(\mathbf{p})] = f^* - \min_{\mathbf{x}_i \in X_i, \forall i \in I} L(\mathbf{p},\mathbf{x}), \quad (2)$$

where the RHS of (2) is the difference in values between the primal optimization problem (1) and the relaxed optimization problem. The CHP method produces a market price $\mathbf{p}^{CHP}$, which minimizes the total uplift payment. In this case, the market price is a solution to the dual problem, $\mathbf{p}^{CHP} \in \arg\max_{\mathbf{p} \in R^T} \min_{\mathbf{x}_i \in X_i, \forall i \in I} L(\mathbf{p},\mathbf{x})$, and the total uplift payment attains its minimum value given by the duality gap [19]-[21]:

$$\min_{\mathbf{p} \in R^T} \sum_{i \in I}[\pi_i^{st.+}(\mathbf{p}) - \pi_i^{st.*}(\mathbf{p})] = \sum_{i \in I}[\pi_i^{st.+}(\mathbf{p}^{CHP}) - \pi_i^{st.*}(\mathbf{p}^{CHP})] = f^* - \max_{\mathbf{p} \in R^T} \min_{\mathbf{x}_i \in X_i, \forall i \in I} L(\mathbf{p},\mathbf{x}).$$

For a given market price $\mathbf{p}$ (not necessarily resulting from the CHP approach), the lost profit of a producer $i$ in the amount of $[\pi_i^{st.+}(\mathbf{p}) - \pi_i^{st.*}(\mathbf{p})]$ is recovered through the uplift payment provided that a producer complies with the centralized dispatch. This payment can be viewed as a modification of the producer's revenue function by adding the term $\delta_{\mathbf{x}_i, \mathbf{x}_i^*}[\pi_i^{st.+}(\mathbf{p}) - \pi_i^{st.*}(\mathbf{p})]$ with a function $\delta_{\mathbf{x}_i, \mathbf{x}_i^*}$ defined as $\delta_{\mathbf{x}_i, \mathbf{x}_i^*} = 1$ for $\mathbf{x}_i = \mathbf{x}_i^*$, and $\delta_{\mathbf{x}_i, \mathbf{x}_i^*} = 0$, otherwise. We have $\delta_{\mathbf{x}_i, \mathbf{x}_i^*} = \delta_{\mathbf{u}_i, \mathbf{u}_i^*} \delta_{\mathbf{g}_i, \mathbf{g}_i^*}$.

Now let us extend the commodity basket of the market and introduce a set $K$ of new products/services described by the inequality constraints $h^k(\mathbf{p},\mathbf{x}) \leq 0$, $k \in K$, with the additively separable functions $h^k(\mathbf{p},\mathbf{x}) = \sum_{i \in I} h_i^k(\mathbf{p},\mathbf{x}_i)$, where $h_i^k(\mathbf{p},\mathbf{x}_i)$ are some functions on $R^T \times X_i$, $\forall k \in K$, $\forall i \in I$. Let us define $h_i(\mathbf{p},\mathbf{x}_i) = \{h_i^1(\mathbf{p},\mathbf{x}_i),...,h_i^{|K|}(\mathbf{p},\mathbf{x}_i)\}$. We allow for the dependence of $h^k(\mathbf{p},\mathbf{x})$ on the market price $\mathbf{p}$, which is just treated as some constant parameter. We also require that the new constraints are redundant on $\Omega$, i.e., $h^k(\mathbf{p},\mathbf{x}) \leq 0$, $\forall x \in \Omega, \forall k \in K$, which means that they do not restrict the feasible set of the original problem. In this case, the centralized dispatch optimization problem is expressed as

$$f^* = \min_{\substack{\mathbf{x}_i \in X_i, \forall i \in I \\ s.t. \\ \sum_{i \in I}\mathbf{g}_i = \mathbf{d} \\ h^k(\mathbf{p},\mathbf{x}) \leq 0, k \in K}} \sum_{i \in I} C_i(\mathbf{x}_i). \quad (3)$$

The Lagrangian relaxation of the new constraints as well as the power balance constraint yields the Lagrange function

$$L(\mathbf{q},\sigma,\mathbf{x}) = \mathbf{q}^T(\mathbf{d} - \sum_{i \in I} \mathbf{g}_i) + \sum_{i \in I} \sigma^T h_i(\mathbf{p},\mathbf{x}_i) + \sum_{i \in I} C_i(\mathbf{x}_i) = \mathbf{q}^T\mathbf{d} - \sum_{i \in I} \pi_i(\mathbf{q},\sigma,\mathbf{x}_i)$$

with multipliers $\mathbf{q} \in R^T$, $\sigma \in R^{|K|}_{\geq 0}$ and a producer $i$ profit function $\pi_i(\mathbf{q},\sigma,\mathbf{x}_i) = \mathbf{q}^T\mathbf{x}_i - \sigma^T h_i(\mathbf{p},\mathbf{x}_i) - C_i(\mathbf{x}_i) = \pi_i^{st.}(\mathbf{q},\mathbf{x}_i) - \sigma^T h_i(\mathbf{p},\mathbf{x}_i)$, $\forall i \in I$. The "vector" of multipliers $\sigma$ is interpreted as the "prices" for the additional products/services $K$. Let us define $\pi_i^*(\mathbf{q},\sigma) = \pi_i(\mathbf{q},\sigma,\mathbf{x}_i^*)$ and $\pi_i^+(\mathbf{q},\sigma) = \max_{\mathbf{x}_i \in X_i} \pi_i(\mathbf{q},\sigma,\mathbf{x}_i)$. For the market price $\mathbf{q} = \mathbf{p}$, this implies

$$f^* - \min_{\mathbf{x}_i \in X_i, \forall i \in I} L(\mathbf{p},\sigma,\mathbf{x}) = \sum_{i \in I}[\pi_i^+(\mathbf{p},\sigma) - \pi_i^*(\mathbf{p},\sigma)] - \sum_{i \in I} \sigma^T h_i(\mathbf{p},\mathbf{x}_i^*). \quad (4)$$

Hence, due to the new term $-\sum_{i \in I} \sigma^T h_i(\mathbf{p},\mathbf{x}_i^*)$ on the RHS of (4), the total uplift payment no longer coincides with the difference between the value of the centralized dispatch optimization problem (3) and the relaxed one. In particular, if the price $\mathbf{p}$ is obtained using the CHP algorithm, then the total uplift payment no longer equals the duality gap [33]. In this case, since the new term is non-negative, the total uplift payment can be lower than the duality gap. We note that the new constraints enter (4) as well as the producer profit functions only through the expression $\sigma^T h_i(\mathbf{p},\mathbf{x}_i)$. Let $\sigma^+$ denote the vector of prices for the commodity set $K$ applied at the market. This implies that the replacement of $|K|$ constraints $h^k(\mathbf{p},\mathbf{x}) \leq 0$ with just one constraint $\sum_{i \in I} \sigma^{+T} h_i(\mathbf{p},\mathbf{x}_i) \leq 0$ with the multiplier equal to 1 results in the same values of both the producer's profit at the centralized dispatch schedule and the producer's maximum profit attainable on its private feasible set as well as the producer's uplift payment for any producer $i$. Consequently, we substitute $|K|$ redundant constraints $h^k(\mathbf{p},\mathbf{x}) \leq 0$ in (3) by just one redundant constraint of the form $\sum_{i \in I} N_i(\mathbf{p},\mathbf{x}_i) \geq 0$ (which holds on $\Omega$) with some functions $N_i(\mathbf{p},\mathbf{x}_i)$ defined on $R^T \times X_i$. Introducing the multiplier $v \geq 0$ to the constraint $\sum_{i \in I} N_i(\mathbf{p},\mathbf{x}_i) \geq 0$, we obtain

$$L(\mathbf{q},v,\mathbf{x}) = \mathbf{q}^T(\mathbf{d} - \sum_{i \in I} \mathbf{g}_i) - v\sum_{i \in I} N_i(\mathbf{p},\mathbf{x}_i) + \sum_{i \in I} C_i(\mathbf{x}_i) = \mathbf{q}^T\mathbf{d} - \sum_{i \in I} \pi_i(\mathbf{q},v,\mathbf{x}_i), \quad (5)$$

with $\pi_i(\mathbf{q},v,\mathbf{x}_i) = \pi_i^{st.}(\mathbf{q},\mathbf{x}_i) + vN_i(\mathbf{p},\mathbf{x}_i)$. Defining $\pi_i^*(\mathbf{q},v) = \pi_i(\mathbf{q},v,\mathbf{x}_i^*)$ and $\pi_i^+(\mathbf{q},v) = \max_{\mathbf{x}_i \in X_i} \pi_i(\mathbf{q},v,\mathbf{x}_i)$ as above, we have for $\mathbf{q} = \mathbf{p}$

$$\sum_{i \in I}[\pi_i^+(\mathbf{p},v) - \pi_i^*(\mathbf{p},v)] = f^* - \min_{\mathbf{x}_i \in X_i, \forall i \in I} L(\mathbf{p},v,\mathbf{x}) - v\sum_{i \in I} N_i(\mathbf{p},\mathbf{x}_i^*). \quad (6)$$

To achieve the total uplift reduction, we need $\sum_{i \in I} N_i(\mathbf{p},\mathbf{x}_i^*) > 0$. Our goal is to minimize the total uplift payment without changing the maximum profit of any producer attainable on its private feasible set. This implies solving the optimization problem

$$\min_{\substack{v \geq 0 \\ s.t. \\ \pi_i^+(\mathbf{p},v) = \pi_i^{st.+}(\mathbf{p}), \forall i \in I}} \sum_{i \in I}[\pi_i^+(\mathbf{p},v) - \pi_i^*(\mathbf{p},v)]. \quad (7)$$

Since $\pi_i^{st.+}(\mathbf{p}) = \pi_i^+(\mathbf{p},0)$, the feasible set of (7) contains $v = 0$ and, therefore, is non-empty. In the next section, we show that replacing the redundancy condition in the form of $\sum_{i \in I} N_i(\mathbf{p},\mathbf{x}_i) \geq 0, \forall \mathbf{x} \in \Omega$, by a stronger condition that implies redundancy of the constraint on a larger set allows solving (7) and explicitly constructing the functions $N_i(\mathbf{p},\mathbf{x}_i)$, $i \in I$, that produce zero total uplift.

## 3. Attaining zero total uplift payment

Consider a set $\times_{i \in I} X_i$, which is a set of $\mathbf{x}$ with $\mathbf{x}_i \in X_i$, $\forall i \in I$. The set $\times_{i \in I} X_i$ would be a feasible set of the problem (1) should the power balance constraint be excluded from the constraint set. Clearly, $\Omega \subset \times_{i \in I} X_i$. Let us replace the redundancy condition $\sum_{i \in I} N_i(\mathbf{p}, \mathbf{x}_i) \geq 0$, $\forall \mathbf{x} \in \Omega$, by a stronger condition

$$\sum_{i \in I} N_i(\mathbf{p}, \mathbf{x}_i) \geq 0, \ \forall \mathbf{x} \in \times_{i \in I} X_i, \qquad (8)$$

which states that the constraint $\sum_{i \in I} N_i(\mathbf{p}, \mathbf{x}_i) \geq 0$ is redundant on $\times_{i \in I} X_i$ not just on $\Omega$.

It is possible to rearrange the terms, which are independent of $\mathbf{x}$, within the sum $\sum_{i \in I} N_i(\mathbf{p}, \mathbf{x}_i)$, so that $\forall i \in I$ we would have $N_i(\mathbf{p}, \mathbf{x}_i) \geq 0$ on $X_i$. Indeed, if this is not the case, define $I_+ = \{i \mid i \in I, \min_{\mathbf{x}_i \in X_i} N_i(\mathbf{p}, \mathbf{x}_i) \geq 0\}$ and $I_- = \{i \mid i \in I, \min_{\mathbf{x}_i \in X_i} N_i(\mathbf{p}, \mathbf{x}_i) < 0\}$. From $\min_{\mathbf{x} \in \times_{i \in I} X_i} \sum_{i \in I} N_i(\mathbf{p}, \mathbf{x}_i) \geq 0$, we obtain $\sum_{i \in I} \min_{\mathbf{x}_i \in X_i} N_i(\mathbf{p}, \mathbf{x}_i) \geq 0$. Introduce

$$\tilde{N}_i(\mathbf{p}, \mathbf{x}_i) = N_i(\mathbf{p}, \mathbf{x}_i) - \min_{\mathbf{x}_i \in X_i} N_i(\mathbf{p}, \mathbf{x}_i), \ \forall i \in I_-,$$

$$\tilde{N}_i(\mathbf{p}, \mathbf{x}_i) = N_i(\mathbf{p}, \mathbf{x}_i) + \frac{\sum_{j \in I_-} \min_{\mathbf{x}_j \in X_j} N_j(\mathbf{p}, \mathbf{x}_j)}{\sum_{k \in I_+} \min_{\mathbf{x}_k \in X_k} N_k(\mathbf{p}, \mathbf{x}_k)} \min_{\mathbf{x}_i \in X_i} N_i(\mathbf{p}, \mathbf{x}_i), \ \forall i \in I_+.$$

Clearly, $\sum_{i \in I} N_i(\mathbf{p}, \mathbf{x}_i) = \sum_{i \in I} \tilde{N}_i(\mathbf{p}, \mathbf{x}_i)$, and $\tilde{N}_i(\mathbf{p}, \mathbf{x}_i) \geq 0$ on $X_i$, $\forall i \in I$. In what follows, we assume that such a rearrangement has been performed resulting in $N_i(\mathbf{p}, \mathbf{x}_i) \geq 0$, $\forall \mathbf{x}_i \in X_i$, $\forall i \in I$. (We note that this condition can be viewed as redundancy of the constraint $N_i(\mathbf{p}, \mathbf{x}_i) \geq 0$ on the private feasible set of a corresponding producer $i$.) For $\forall i \in I$, $\forall v \geq 0$, this yields $\pi_i(\mathbf{p}, v, \mathbf{x}_i) \geq \pi_i^{st.}(\mathbf{p}, \mathbf{x}_i)$, $\forall \mathbf{x}_i \in X_i$; therefore, $\pi_i^+(\mathbf{p}, v) \geq \pi_i^{st.+}(\mathbf{p})$. Now we show that the feasible set of (7) is given by a set of optimal points of a certain dual problem. Because the constraint $\sum_{i \in I} N_i(\mathbf{p}, \mathbf{x}_i) \geq 0$ is redundant on $\times_{i \in I} X_i$, the problem

$$\min_{\mathbf{x}_i \in X_i, \forall i \in I} L(\mathbf{p}, v, \mathbf{x}) = \mathbf{p}^T \mathbf{d} - \sum_{i \in I} \max_{\mathbf{x}_i \in X_i} [\pi_i^{st.}(\mathbf{p}, \mathbf{x}_i) + v N_i(\mathbf{p}, \mathbf{x}_i)]$$

can be viewed as resulting from the problem $\min_{\substack{\mathbf{x}_i \in X_i, \forall i \in I \\ \sum_{i \in I} N_i(\mathbf{p}, \mathbf{x}_i) \geq 0}} L(\mathbf{p}, 0, \mathbf{x})$ after the Lagrangian relaxation of this redundant constraint. From $\min_{\substack{\mathbf{x}_i \in X_i, \forall i \in I \\ \sum_{i \in I} N_i(\mathbf{p}, \mathbf{x}_i) \geq 0}} L(\mathbf{p}, 0, \mathbf{x}) = \min_{\mathbf{x}_i \in X_i, \forall i \in I} L(\mathbf{p}, 0, \mathbf{x})$ and $\min_{\mathbf{x}_i \in X_i, \forall i \in I} L(\mathbf{p}, 0, \mathbf{x}) \leq \max_{v \geq 0} \min_{\mathbf{x}_i \in X_i, \forall i \in I} L(\mathbf{p}, v, \mathbf{x}) \leq \min_{\mathbf{x}_i \in X_i, \forall i \in I} L(\mathbf{p}, 0, \mathbf{x})$, it follows that

$$\max_{v \geq 0} \min_{\mathbf{x}_i \in X_i, \forall i \in I} L(\mathbf{p}, v, \mathbf{x}) = \min_{\mathbf{x}_i \in X_i, \forall i \in I} L(\mathbf{p}, 0, \mathbf{x}), \qquad (9)$$

which implies $\min_{v \geq 0} \sum_{i \in I} \pi_i^+(\mathbf{p}, v) = \sum_{i \in I} \pi_i^{st.+}(\mathbf{p})$. Let $v^+$ be an optimal solution to (9), then from $\pi_i^+(\mathbf{p}, v) \geq \pi_i^{st.+}(\mathbf{p})$, $\forall v \geq 0$, we obtain $\pi_i^+(\mathbf{p}, v^+) = \pi_i^{st.+}(\mathbf{p})$, $\forall i \in I$.

This has two implications. First, the introduction and the Lagrangian relaxation of the constraint that is redundant on $\times_{i \in I} X_i$ introduce no duality gap and, as a result, do not change the duality gap originating from the power balance constraint dualization. Second, the feasible set of (7) is given by the set of maximizers of the dual problem present on the LHS of (9). Since $\min_{\mathbf{x}_i \in X_i, \forall i \in I} L(\mathbf{p}, v, \mathbf{x})$ is a concave function of $v$ as a pointwise minimum of the concave (linear) functions, the set of maximizers is convex. If $\sum_{i \in I} N_i(\mathbf{p}, \mathbf{x}_i^*) > 0$, then this set is bounded (otherwise, one could achieve arbitrary large $\pi_i^+(\mathbf{p}, v)$ for $i$ with $N_i(\mathbf{p}, \mathbf{x}_i^*) > 0$, which is incompatible with $\pi_i^+(\mathbf{p}, v^+) = \pi_i^{st.+}(\mathbf{p})$). Since the function $\min_{\mathbf{x}_i \in X_i, \forall i \in I} L(\mathbf{p}, v, \mathbf{x})$ is continuous on $v \in [0, +\infty)$, its superlevel sets are closed. As a result, if $\sum_{i \in I} N_i(\mathbf{p}, \mathbf{x}_i^*) > 0$, the set of maximizers has the form $[0, v_{\max}^+(\mathbf{p})]$ with some $v_{\max}^+(\mathbf{p}) \geq 0$ and from (6) and (9) we obtain

$$\min_{\substack{v \geq 0 \\ s.t. \\ \pi_i^+(\mathbf{p},v)=\pi_i^{st.+}(\mathbf{p}), \forall i \in I}} \sum_{i \in I}[\pi_i^+(\mathbf{p},v) - \pi_i^*(\mathbf{p},v)] = f^* - \min_{\mathbf{x}_i \in X_i, \forall i \in I} L(\mathbf{p},0,\mathbf{x}) - v_{\max}^+(\mathbf{p})\sum_{i \in I} N_i(\mathbf{p},\mathbf{x}_i^*).$$

Thus, if $v_{\max}^+(\mathbf{p}) = 0$, there is no total uplift reduction. If $v_{\max}^+(\mathbf{p}) > 0$, then $v_{\max}^+(\mathbf{p}) \neq 1$ can be always absorbed in $\sum_{i \in I} N_i(\mathbf{p},\mathbf{x}_i^*)$ by rescaling the functions $N_i(\mathbf{p},\mathbf{x}_i)$, $i \in I$. Therefore, we may focus on the functions $N_i(\mathbf{p},\mathbf{x}_i)$, $i \in I$, that yield $v_{\max}^+(\mathbf{p}) = 1$. We have the following proposition that formulates sufficient conditions for the functions $N_i(\mathbf{p},\mathbf{x}_i)$, $i \in I$, with $v_{\max}^+(\mathbf{p}) = 1$ to eliminate the total uplift payment.

**Proposition** Let $\forall i \in I$ the function $N_i(\mathbf{p},\mathbf{x}_i)$ satisfy the following conditions:

$$\pi_i^{st.+}(\mathbf{p}) = \max_{\mathbf{x}_i \in X_i}[\pi_i^{st.}(\mathbf{p},\mathbf{x}_i) + N_i(\mathbf{p},\mathbf{x}_i)], \quad (10)$$

$$\pi_i^{st.+}(\mathbf{p}) = \pi_i^{st.}(\mathbf{p},\mathbf{x}_i^*) + N_i(\mathbf{p},\mathbf{x}_i^*), \quad (11)$$

$$N_i(\mathbf{p},\mathbf{x}_i) \geq 0, \ \forall \mathbf{x}_i \in X_i, \quad (12)$$

then

- the constraint $\sum_{i \in I} N_i(\mathbf{p},\mathbf{x}_i) \geq 0$ is redundant on $\Omega$;
- $v = 1$ is the optimal point of (7) with $\sum_{i \in I}[\pi_i^+(\mathbf{p},1) - \pi_i^*(\mathbf{p},1)] = 0$, i.e., the total uplift payment vanishes.

*Proof.* The statement of first bullet readily follows from the redundancy of $\sum_{i \in I} N_i(\mathbf{p},\mathbf{x}_i) \geq 0$ on $\times_{i \in I} X_i$ due to (12). To prove the claim of the second bullet, we note that since $\sum_{i \in I}[\pi_i^+(\mathbf{p},v) - \pi_i^*(\mathbf{p},v)] \geq 0$, $\forall v \geq 0$, any feasible $v$ with $\sum_{i \in I}[\pi_i^+(\mathbf{p},v) - \pi_i^*(\mathbf{p},v)] = 0$ is an optimal point of (7). Thus, it suffices to show that $v = 1$ is feasible in (7) and gives zero total uplift payment. Feasibility of $v = 1$ follows from (10), which implies $\pi_i^{st.+}(\mathbf{p}) = \pi_i^+(\mathbf{p},1)$, $\forall i \in I$, while $\sum_{i \in I}[\pi_i^+(\mathbf{p},1) - \pi_i^*(\mathbf{p},1)] = 0$ is entailed from (10) and (11). Proposition is proved.

The conditions (10) - (12) imply that a non-negative function $N_i(\mathbf{p},\mathbf{x}_i)$ vanishes on $\arg\max_{\mathbf{x}_i \in X_i} \pi_i^{st.}(\mathbf{p},\mathbf{x}_i)$, at $\mathbf{x}_i = \mathbf{x}_i^*$ it equals the producer's original uplift payment $\pi_i^{st.+}(\mathbf{p}) - \pi_i^{st.}(\mathbf{p},\mathbf{x}_i^*)$, and it has sufficiently low values on the rest of $X_i$ to ensure that the profit function $\pi_i^{st.}(\mathbf{p},\mathbf{x}_i) + N_i(\mathbf{p},\mathbf{x}_i)$ is bounded by $\pi_i^{st.+}(\mathbf{p})$. The general form expression of such a function is given by

$$N_i(\mathbf{p},\mathbf{x}_i) = \min[\pi_i^{st.+}(\mathbf{p}) - \pi_i^{st.}(\mathbf{p},\mathbf{x}_i); \delta_{\mathbf{x}_i, \mathbf{x}_i^*}[\pi_i^{st.+}(\mathbf{p}) - \pi_i^{st.*}(\mathbf{p})] + \gamma_i(\mathbf{p},\mathbf{x}_i)], \ \forall \mathbf{x}_i \in X_i, \quad (13)$$

with an arbitrary non-negative function $\gamma_i(\mathbf{p},\mathbf{x}_i)$: $R^T \times X_i \to R$, $\gamma_i(\mathbf{p},\mathbf{x}_i) \geq 0$, $\forall \mathbf{x}_i \in X_i$. In [34] it was shown that a function $N_i(\mathbf{p},\mathbf{x}_i)$ satisfies (10) – (12) if and only if it has the form (13). Thus, the freedom to choose an arbitrary non-negative function $\gamma_i(\mathbf{p},\mathbf{x}_i)$ for each $N_i(\mathbf{p},\mathbf{x}_i)$ entails that there could be the multiple sets of functions $N_i(\mathbf{p},\mathbf{x}_i)$, $i \in I$, producing zero total uplift payment without affecting the maximum value of the profit function of any producer. This translates into the multiple choices for the redundant constraint $\sum_{i \in I} N_i(\mathbf{p},\mathbf{x}_i) \geq 0$.

We also note that neither in Proposition nor in equivalence of (10) – (12) and (13) we used the fact that the revenue function $R_i^{st.}(\mathbf{p},\mathbf{x}_i)$, which enters the profit function $\pi_i^{st.}(\mathbf{p},\mathbf{x}_i) = R_i^{st.}(\mathbf{p},\mathbf{x}_i) - C_i(\mathbf{x}_i)$, had the form of $R_i^{st.}(\mathbf{p},\mathbf{x}_i) = \mathbf{p}^T \mathbf{g}_i$. Therefore, the expression (13) is also applicable in the general case of possibly non-linear revenue functions $R_i^{st.}(\mathbf{p},\mathbf{x}_i)$, $i \in I$, with $\mathbf{p}$ being a set of the pricing parameters. In these cases, each redundant constraint $N_i(\mathbf{p},\mathbf{x}_i) \geq 0$ is introduced in the profit optimization problem of a producer $i$, and we also have the total uplift optimization problem (7) with $\pi_i(\mathbf{p},v,\mathbf{x}_i) = \pi_i^{st.}(\mathbf{p},\mathbf{x}_i) + vN_i(\mathbf{p},\mathbf{x}_i)$, $\pi_i^*(\mathbf{p},v) = \pi_i(\mathbf{p},v,\mathbf{x}_i^*)$, and $\pi_i^+(\mathbf{p},v) = \max_{\mathbf{x}_i \in X_i} \pi_i(\mathbf{p},v,\mathbf{x}_i)$ with some functions $N_i(\mathbf{p},\mathbf{x}_i)$, $i \in I$. If the latter are given by (13), then the conditions (10) – (12) are satisfied and the statement of Proposition holds ensuring zero total uplift payment at $v = 1$.

## 4. Applications

In this section, we apply the findings to the cases of the uninode multi-period power markets with marginal pricing and CHP. We also show how the suggested approach can be applied to a multi-node power system with FTR holders (which is illustrated by the example of a two-node power system).

### 4.1. Power markets with marginal pricing

In the case of a power market with marginal pricing, $\mathbf{p}$ is identified as the marginal price vector faced by a producer $i$ with the revenue function $R_i^{st.}(\mathbf{p},\mathbf{x}_i) = \mathbf{p}^T \mathbf{g}_i$. Let us assume that for each fixed value of $\mathbf{u}_i$ both the generator private feasible set and the cost function are convex. The marginal pricing method entails

$$\pi_i^{st.*}(\mathbf{p}) = \max_{\substack{\mathbf{x}_i \in X_i \\ s.t.\ \mathbf{u}_i = \mathbf{u}_i^*}} \pi_i^{st.}(\mathbf{p},\mathbf{x}_i). \quad (14)$$

As it was mentioned above, the choice for $N_i(\mathbf{p},\mathbf{x}_i)$ satisfying (10) – (12) may not be unique. Below we provide an expression for $N_i(\mathbf{p},\mathbf{x}_i)$, which originates from the constraint $1 - u_i^t \geq 0$ if $u_i^{*t} = 0$ and the constraint $u_i^t \geq 0$ if $u_i^{*t} = 1$. Such a choice can be written as $(u_i^t)^{u_i^{*t}}(1-u_i^t)^{(1-u_i^{*t})} \geq 0$, which is redundant on the private feasible set of a producer $i$. Cleary, we have $(u_i^t)^{u_i^{*t}}(1-u_i^t)^{(1-u_i^{*t})} = \delta_{\mathbf{u}_i,\mathbf{u}_i^*}$. It is straightforward to check that the multiplier associated with the redundant constraint $\delta_{\mathbf{u}_i,\mathbf{u}_i^*} \geq 0$ in a producer profit optimization problem is given by $\pi_i^{st.+}(\mathbf{p}) - \pi_i^{st.*}(\mathbf{p})$. Absorbing the multiplier in the constraint, we arrive at $N_i(\mathbf{p},\mathbf{x}_i) = \delta_{\mathbf{u}_i,\mathbf{u}_i^*}[\pi_i^{st.+}(\mathbf{p}) - \pi_i^{st.*}(\mathbf{p})]$, which can be obtained from (13) with $\gamma_i(\mathbf{p},\mathbf{x}_i) = (\delta_{\mathbf{u}_i,\mathbf{u}_i^*} - \delta_{\mathbf{x}_i,\mathbf{x}_i^*})[\pi_i^{st.+}(\mathbf{p}) - \pi_i^{st.*}(\mathbf{p})]$. It can be verified that the conditions (10) – (12) hold. Indeed, the function $N_i(\mathbf{p},\mathbf{x}_i)$ is non-negative, at $\mathbf{x}_i = \mathbf{x}_i^*$ it has the value of the producer's $i$ original uplift, and it does not affect the maximum value of the producer profit function:

$$\max_{\mathbf{x}_i \in X_i}[\pi_i^{st.}(\mathbf{p},\mathbf{x}_i) + N_i(\mathbf{p},\mathbf{x}_i)] = \max\{\max_{\substack{\mathbf{x}_i \in X_i \\ s.t.\ \mathbf{u}_i \neq \mathbf{u}_i^*}} \pi_i^{st.}(\mathbf{p},\mathbf{x}_i); \pi_i^{st.+}(\mathbf{p}) - \pi_i^{st.*}(\mathbf{p}) + \max_{\substack{\mathbf{x}_i \in X_i \\ s.t.\ \mathbf{u}_i = \mathbf{u}_i^*}} \pi_i^{st.}(\mathbf{p},\mathbf{x}_i)\} = \pi_i^{st.+}(\mathbf{p}).$$

The redundant constraint for the centralized dispatch optimization problem is given by $\sum_{i \in I} N_i(\mathbf{p},\mathbf{x}_i) = \sum_{i \in I} \delta_{\mathbf{u}_i,\mathbf{u}_i^*}[\pi_i^{st.+}(\mathbf{p}) - \pi_i^{st.*}(\mathbf{p})] \geq 0$. The resulting profit function of a producer $i$ is expressed as $\pi_i(\mathbf{p},\mathbf{x}_i) = \mathbf{p}^T \mathbf{g}_i - C_i(\mathbf{x}_i) + \delta_{\mathbf{u}_i,\mathbf{u}_i^*}[\pi_i^{st.+}(\mathbf{p}) - \pi_i^{st.*}(\mathbf{p})]$. We note that for the given market price $\mathbf{p}$ and statuses of the unit, the new term in the profit function is constant.

*Example 1.* Consider a uninode single-period market with marginal pricing. The demand is assumed to be fixed. A private feasible set of a producer $i$ is given by $X_i = \{(u_i, g_i) \mid u_i \in \{0,1\}, g_i \in R, u_i g_i^{\min} \leq g_i \leq u_i g_i^{\max}\}$. For each producer $i \in I$, we choose $\gamma_i(p,x_i) = (\delta_{u_i,u_i^*} - \delta_{x_i,x_i^*})[\pi_i^{st.+}(p) - \pi_i^{st.*}(p)]$ as above. Thus, $N_i(p,x_i) = \delta_{u_i,u_i^*}[\pi_i^{st.+}(p) - \pi_i^{st.*}(p)]$, $\forall i \in I$. Clearly, $N_i(p,x_i) \geq 0$, $\forall x_i \in X_i$, $\forall i \in I$, and the constraint $N_i(p,x_i) \geq 0$ is redundant on the private feasible set of producer $i$. We have $\delta_{u_i,0} = (1-u_i)$, $\delta_{u_i,1} = u_i$. Consequently, for $u_i^* = 0$ the constraint $N_i(p,x_i) \geq 0$ up to a non-negative factor coincides with the redundant constraint $u_i \leq 1$, while for $u_i^* = 1$ the constraint $N_i(p,x_i) \geq 0$ up to a non-negative factor coincides with the redundant constraint $u_i \geq 0$. If the generator has the lost profit (i.e., it is offline in the centralized dispatch solution with $\pi_i^{st.*}(p) = 0$, but $\pi_i^{st.+}(p) > 0$), then $N_i(p, x_i) = (1-u_i)\pi_i^{st.+}(p)$, which compensates the producer for the lost profit $\pi_i^{st.+}(p)$ if it complies with the offline status set by the centralized dispatch solution. Likewise, if the generator operates at a loss (i.e., the generator is online in the centralized dispatch solution, but $\pi_i^{st.*}(p) < 0$ so that it is not recovering its cost at the given market price $p$), then $\pi_i^{st.+}(p) = 0$ and $N_i(p, x_i) = -u_i \pi_i^{st.*}(p)$, which implies compensating the producer if it stays online. Due to (14), the maximum profit of the producer for a given online status of its generating unit is attained at the centralized dispatch schedule $x_i^*$. Thus, in this case, the producer has economic incentives to deliver $x_i^*$, which results in zero producer's profit implying full recovery of its cost. As a result, $\forall i \in I$ the function $N_i(p, x_i)$ satisfies (10) – (12); therefore, the constraint $\sum_{i \in I} N_i(p, x_i) \geq 0$ is redundant on the feasible set of (1) and results in zero total uplift.

### 4.2. Power markets with CHP pricing

In the case of the CHP method, both the power balance constraint and the redundant constraint are dualized to obtain the corresponding prices. Due to the redundant nature of the new constraint, the duality gap is either lowered or remains the same [32]. However, it was shown [32] that dualization of the linear redundant constraint does not affect the duality gap. This motivated the study [33], where it was shown that the total uplift payment is reduced (but is generally non-zero) if a linear redundant constraint (which is a linear combination of the power balance constraint and the capacity limit constraints of all the producers) for each time instance is introduced in the centralized dispatch optimization problem. It is possible to illustrate that in the general case the linear redundant constraints are not able to eliminate the total uplift payment.

Our findings show that it is possible to substitute the linearity requirement by a redundancy condition (8) that also ensures invariance of the duality gap if the (possibly non-linear) redundant constraint is dualized. Thus, the resulting duality gap originates from the power balance constraint dualization and is not affected by the new constraint. Moreover, the set of market prices for power is unaltered as well. Indeed, from (9) for $\mathbf{p} = \mathbf{p}^{CHP}$ it follows that if $(\mathbf{p}^{CHP}, v^+)$ is a dual problem optimal point, then so is $(\mathbf{p}^{CHP}, 0)$. In Section 3, we have demonstrated that the condition (8) allows sufficient freedom to fully eliminate the total uplift payment without changing the duality gap.

The uplift payment $\pi_i^{st.+}(\mathbf{p}^{CHP}) - \pi_i^{st.*}(\mathbf{p}^{CHP})$ can be viewed as the cost of the commitment ticket payable to the generator $i$ for following the centralized dispatch [19]-[21]. This can be interpreted as modifying a producer $i$ revenue function to include the term $\delta_{\mathbf{x}_i, \mathbf{x}_i^*}[\pi_i^{st.+}(\mathbf{p}^{CHP}) - \pi_i^{st.*}(\mathbf{p}^{CHP})]$, which resembles the term appearing in the profit function after the Lagrangian relaxation of the constraint $\delta_{\mathbf{x}_i, \mathbf{x}_i^*} \geq 0$, which is redundant on $X_i$, with the associated multiplier $[\pi_i^{st.+}(\mathbf{p}^{CHP}) - \pi_i^{st.*}(\mathbf{p}^{CHP})]$. Indeed, setting $\gamma_i(\mathbf{p}^{CHP}, \mathbf{x}_i) \equiv 0$ in (13) yields $N_i(\mathbf{p}^{CHP}, \mathbf{x}_i) = \delta_{\mathbf{x}_i, \mathbf{x}_i^*}[\pi_i^{st.+}(\mathbf{p}^{CHP}) - \pi_i^{st.*}(\mathbf{p}^{CHP})]$.

Consequently, the original total uplift payment $\sum_{i \in I} \delta_{\mathbf{x}_i, \mathbf{x}_i^*}[\pi_i^{st.+}(\mathbf{p}^{CHP}) - \pi_i^{st.*}(\mathbf{p}^{CHP})]$ can be expressed as the term originating from the Lagrangian relaxation of the redundant constraint $\sum_{i \in I} N_i(\mathbf{p}^{CHP}, \mathbf{x}_i) \geq 0$ in the centralized dispatch optimization problem. Therefore, the introduction and the Lagrangian relaxation of this constraint fully eliminate the total uplift payment if a market price is set by the CHP algorithm. The alternative choices for the functions $\gamma_i(\mathbf{p}^{CHP}, \mathbf{x}_i)$ for the markets with the price set by the CHP method were given in [34].

### 4.3 Two-node power markets with FTR-holders

Our approach can be readily extended to multi-node power systems. The Lagrangian relaxation of the power flow constraints leads to the terms attributable to the FTR holders. As an example, consider a multi-period two-node power system with the vector of power flows denoted as $\mathbf{F}$, $\mathbf{F} \in R^T$. (We adopt the convention that positive flow is in the direction to node #2.) Let us define the sets $I_1$, $I_2$, which include the producers located in node #1 and node #2, respectively, the vectors of nodal demand $\mathbf{d}_1$ and $\mathbf{d}_2$, the nodal market prices $\mathbf{p}_1$ and $\mathbf{p}_2$, and the vector of maximum power flows $\mathbf{F}^{\max}$. In this case, the centralized dispatch optimization problem has the form

$$f^* = \min_{\substack{\mathbf{x}_i \in X_i, \forall i \in I \\ \mathbf{F} \in R^T \\ s.t. \\ -\mathbf{F}^{\max} \leq \mathbf{F} \leq \mathbf{F}^{\max} \\ \sum_{i \in I_1} \mathbf{g}_i = \mathbf{d}_1 + \mathbf{F}, \\ \sum_{i \in I_2} \mathbf{g}_i + \mathbf{F} = \mathbf{d}_2}} \sum_{i \in I} C_i(\mathbf{x}_i) \cdot \quad (15)$$

with the optimal power flow denoted as $\mathbf{F}^*$. Let us introduce the redundant constraint of the form $\sum_{i \in I_1} N_i(\mathbf{p}_1, \mathbf{x}_i) + \sum_{i \in I_2} N_i(\mathbf{p}_2, \mathbf{x}_i) + N_{FTR(12)}(\mathbf{p}_1, \mathbf{p}_2, \mathbf{F}) \geq 0$ with $N_i \geq 0$, $\forall \mathbf{x}_i \in X_i$, $\forall i \in I$, and $N_{FTR(12)}(\mathbf{p}_1, \mathbf{p}_2, \mathbf{F}) \geq 0$, $\forall \mathbf{F}: -\mathbf{F}^{\max} \leq \mathbf{F} \leq \mathbf{F}^{\max}$. The Lagrangian relaxation of this constraint as well as the power balance constraints gives

$$\min_{\substack{\mathbf{x}_i \in X_i, \forall i \in I \\ -\mathbf{F}^{\max} \leq \mathbf{F} \leq \mathbf{F}^{\max}}} L(\mathbf{q}_1, \mathbf{q}_2, \nu, \mathbf{x}, \mathbf{F}) = \mathbf{q}_1^T \mathbf{d}_1 + \mathbf{q}_2^T \mathbf{d}_2 - \sum_{i \in I_1} \max_{\mathbf{x}_i \in X_i} [\pi_i^{st.}(\mathbf{q}_1, \mathbf{x}_i) + \nu N_i(\mathbf{p}_1, \mathbf{x}_i)] - \sum_{i \in I_2} \max_{\mathbf{x}_i \in X_i} [\pi_i^{st.}(\mathbf{q}_2, \mathbf{x}_i) + \nu N_i(\mathbf{p}_2, \mathbf{x}_i)]$$
$$- \max_{-\mathbf{F}^{\max} \leq \mathbf{F} \leq \mathbf{F}^{\max}} [\pi_{FTR(12)}^{st.}(\mathbf{q}_1, \mathbf{q}_2, \mathbf{F}) + \nu N_{FTR(12)}(\mathbf{p}_1, \mathbf{p}_2, \mathbf{F})]$$

where $\pi_{FTR(12)}^{st.}(\mathbf{q}_1, \mathbf{q}_2, \mathbf{F}) = (\mathbf{q}_2 - \mathbf{q}_1)^T \mathbf{F}$ is profit of the FTR holders. Let us define the following functions associated with the FTR holders $\pi_{FTR(12)}(\mathbf{q}_1, \mathbf{q}_2, \nu, \mathbf{F}) = \pi_{FTR(12)}^{st.}(\mathbf{q}_1, \mathbf{q}_2, \mathbf{F}) + \nu N_{FTR(12)}(\mathbf{p}_1, \mathbf{p}_2, \mathbf{F})$, $\pi_{FTR(12)}^*(\mathbf{q}_1, \mathbf{q}_2, \nu) = \pi_{FTR(12)}(\mathbf{q}_1, \mathbf{q}_2, \nu, \mathbf{F}^*)$, $\pi_{FTR(12)}^+(\mathbf{q}_1, \mathbf{q}_2, \nu) = \max_{-\mathbf{F}^{\max} \leq \mathbf{F} \leq \mathbf{F}^{\max}} \pi_{FTR(12)}(\mathbf{q}_1, \mathbf{q}_2, \nu, \mathbf{F})$, and $\pi_{FTR(12)}^{st.+}(\mathbf{q}_1, \mathbf{q}_2) = \max_{-\mathbf{F}^{\max} \leq \mathbf{F} \leq \mathbf{F}^{\max}} \pi_{FTR(12)}^{st.}(\mathbf{q}_1, \mathbf{q}_2, \mathbf{F})$. For $\mathbf{q}_1 = \mathbf{p}_1$, $\mathbf{q}_2 = \mathbf{p}_2$, the total uplift payment reads

$$\sum_{i \in I_1} [\pi_i^+(\mathbf{p}_1, \nu) - \pi_i^*(\mathbf{p}_1, \nu)] + \sum_{i \in I_2} [\pi_i^+(\mathbf{p}_2, \nu) - \pi_i^*(\mathbf{p}_2, \nu)] + \pi_{FTR(12)}^+(\mathbf{p}_1, \mathbf{p}_2, \nu) - \pi_{FTR(12)}^*(\mathbf{p}_1, \mathbf{p}_2, \nu)$$
$$= f^* - \min_{\substack{\mathbf{x}_i \in X_i, \forall i \in I \\ -\mathbf{F}^{\max} \leq \mathbf{F} \leq \mathbf{F}^{\max}}} L(\mathbf{p}_1, \mathbf{p}_2, \nu, \mathbf{x}, \mathbf{F}) - \nu[\sum_{i \in I_1} N_i(\mathbf{p}_1, \mathbf{x}_i^*) + \sum_{i \in I_2} N_i(\mathbf{p}_2, \mathbf{x}_i^*) + N_{FTR(12)}(\mathbf{p}_1, \mathbf{p}_2, \mathbf{F}^*)] \quad (16)$$

To account for the FTR holders, the constraint set of (7) is extended to include the constraint $\pi_{FTR(12)}^{st.+}(\mathbf{p}_1, \mathbf{p}_2) = \pi_{FTR(12)}^+(\mathbf{p}_1, \mathbf{p}_2, \nu)$. If (10) – (12) hold for both the producers and FTR holders, then (16) for $\nu = 1$ reads

$$f^* - \min_{\substack{\mathbf{x}_i \in X_i, \forall i \in I \\ -\mathbf{F}^{\max} \leq \mathbf{F} \leq \mathbf{F}^{\max}}} L(\mathbf{p}_1, \mathbf{p}_2, 0, \mathbf{x}, \mathbf{F}) - \sum_{i \in I_1} N_i(\mathbf{p}_1, \mathbf{x}_i^*) - \sum_{i \in I_2} N_i(\mathbf{p}_2, \mathbf{x}_i^*) - N_{FTR(12)}(\mathbf{p}_1, \mathbf{p}_2, \mathbf{F}^*)$$

and vanishes due to the statement of Proposition.

*Example 2.* Consider a two-node single-period power market with the fixed load. Let the system have two producers (one in each node) with the linear cost functions of the form $C_i(x_i) = a_i g_i + w_i u_i$ defined on the private feasible sets $X_i = \{(u_i, g_i) \mid u_i \in \{0,1\}, g_i \in R, u_i g_i^{\min} \leq g_i \leq u_i g_i^{\max}\}$, $i \in \{1,2\}$. The fixed load $d$ of 150 MWh is assigned to node #1. The feasible set for the power flow variable $F$ (with the direction from node 1 to node 2 being positive) is given by $-F^{\max} \leq F \leq F^{\max}$ with $F^{\max} = 50$ MWh, which corresponds to the transmission capacity of 50 MW. The parameters of the model are given in Table 1.

**Table 1. Parameter values for Example 2**

|  | $g_i^{\min}$, MWh | $g_i^{\max}$, MWh | $a_i$, \$/MWh | $w_i$, \$ |
|---|---|---|---|---|
| Producer 1 | 100 | 200 | 15.00 | 20.00 |
| Producer 2 | 150 | 200 | 10.00 | 0.00 |

The power balance constraints are expressed as $g_1 = d + F$, $g_2 + F = 0$. The centralized dispatch solution is given by $u_1^* = 1$, $g_1^* = 150 MWh$, $u_2^* = g_2^* = F^* = 0$. Let us set the market prices using the CHP method, which gives the nodal market prices $p_1 = \$15.10/MWh$ and $p_2 = \$10.00/MWh$. Thus, $\pi_1^{st.+}(p_1) = 0$, $\pi_1^{st.*}(p_1) = -\$5$, $\pi_2^{st.+}(p_2) = \pi_2^{st.*}(p_2) = 0$, $\pi_{FTR(12)}^{st.+}(p_1, p_2) = \$255$, $\pi_{FTR(12)}^{st.*}(p_1, p_2) = 0$. As a result, the uplift payment for producer 2 vanishes, while the uplift payments for producer 1 and the FTR holders equal \$5 and \$255, respectively. For producer 1, let us choose the function $\gamma_1(p_1, x_1)$ in (13) so that the function $\delta_{x_1, x_1^*}[\pi_1^{st.+}(p_1) - \pi_1^{st.*}(p_1)] + \gamma_1(p_1, x_1)$ has no discontinuity. This can be achieved, for example, by setting $\delta_{x_1, x_1^*}[\pi_1^{st.+}(p_1) - \pi_1^{st.*}(p_1)] + \gamma_1(p_1, x_1) = \min[5, (200 - g_1)/10]$, which entails

$$\gamma_1(p_1, x_1) = 5\{\min[1, (200 - g_1)/50] - \delta_{150, g_1}\}, \quad N_1(p_1, x_1) = 0.1 \min[200u_1 - g_1; 50; 200 - g_1].$$

Likewise, let us choose $\gamma_{FTR(12)}(p_1, p_2, F) = 255 \min[F/50 + 1 - \delta_{F,0}; 1 - \delta_{F,0}]$, which gives $N_{FTR(12)}(p_1, p_2, F) = 255 \min[1 + F/50; 1]$. In this case, just one redundant constraint of the form $N_1(p_1, x_1) + N_{FTR(12)}(p_1, p_2, F) \geq 0$ suffices to eliminate the total uplift payment. The consideration can be straightforwardly generalized to a multi-node power market.

## 5. Conclusion

Integration of distributed energy resources to a power market is one of the most notable changes over the last decade. Utilization of the distributed renewable energy sources leads to the operational challenges related to both ramping and reserve needs of power systems due to an increased volatility of net demand and the cost recovery issues related to low or even negative market prices for power. Combined with zero marginal cost of output for the distributed renewable generation, these factors make significant contributions to the total lost profit of the dispatchable generation that needs to be compensated to ensure the economic stability of the centralized dispatch outcome.

In the framework of a uninode multi-period power market with fixed demand, we studied the influence of the Lagrangian relaxation of the additively separable redundant constraints (i.e., the additively separable constraints that hold on the feasible set of the centralized dispatch optimization problem) on the total uplift payment, which compensates the total lost profit of the market players. We showed that a set of such constraints can be replaced with just one redundant constraint having the same effect on the total uplift payment. In the general case, such a procedure may change the (post-uplift) profit of the market players. The previous research on the subject demonstrated that the linear redundant constraint may lower the total uplift payment without affecting the (post-uplift) profit received by the market players. However, the linearity condition is generally too restrictive to attain zero total uplift. This leads to the question if it is possible to construct a (possibly non-linear) redundant constraint that reduces the total uplift payment to zero (which, in turn, implies vanishing uplift payment to each market player) provided that the (post-uplift) profit received by each market player is unchanged by the Lagrangian relaxation of the new constraint.

For this purpose, we considered a special type of redundant constraints – the constraints that hold not only on the feasible set of the centralized dispatch optimization problem but also on the larger set obtained by relaxing the power balance constraint. The Lagrangian relaxation of such a constraint does not affect the market player profits provided that the associated multiplier belongs to a certain set. In particular, the duality gap, which originates from the power balance constraint dualization, is unaltered.

Imposing the condition of zero total uplift payment, we arrived at the general expression for just one additively separable (generally non-linear) constraint that belongs to this class of the redundant constraints and eliminates the total uplift payment after the Lagrangian relaxation procedure. We explicitly constructed a family of such redundant constraints parameterized by a set of arbitrary non-negative functions. If the uniform market price for power is set by the CHP method, the dualization of both the power balance constraint and the new redundant constraint has no effect on the duality gap, results in the same set of the market prices and the same maximum profit for each market player but gives zero total uplift payment.

The proposed approach is rather general and is applicable to the cases with non-uniform or non-linear pricing as well. (In these cases, the new redundant constraints are introduced in the private profit optimization problems of the market players.)

The analysis can be straightforwardly extended to multi-node multi-period markets with price-sensitive demand.